\documentclass[12pt,a4paper]{amsart}

\usepackage[margin=1in,includeheadfoot]{geometry}

\usepackage{tikz}
\tikzset{every picture/.style = {inner sep=0pt,baseline}}
\tikzset{p/.style = {draw, shape          = circle,
                                 text           = black,
                                 font=\tiny,
                                 inner sep      = 0pt,
                                 outer sep      = 0pt,
                                 minimum size   = 1pt}}

\tikzset{r/.style = {red, very thick}}
\tikzset{g/.style = {green, very thick}}
\tikzset{b/.style = {blue, very thick}}
\tikzset{o/.style = {orange, very thin}}
\tikzset{c/.style = {cyan, very thin}}

\newtheorem{df}{Definition}[section]
\newtheorem*{thm*}{Theorem}
\newtheorem{thm}{Theorem}[section]

\newtheorem{cly}{Corollary}[section]

\usepackage{euscript}

\newcommand{\C}{\mathbb{C}}
\newcommand{\R}{\mathbb{R}}

\newcommand{\Z}{\mathbb{Z}}
\newcommand{\frakg}{\mathfrak{g}}
\newcommand{\frakt}{\mathfrak{t}}
\newcommand{\frakit}{\mathfrak{it}}
\newcommand{\frakD}{\mathfrak{D}}
\newcommand{\Ad}{\mathrm{Ad}}
\newcommand{\Tr}{\mathrm{Tr}}

\begin{document}
\title[Maximal order Abelian subgroups]
{Maximal order Abelian subgroups of Coxeter groups.}

\author{John M. Burns and Goetz Pfeiffer}
 \maketitle
{ \bf Abstract:} In this note we give a classification of the
Maximal order Abelian subgroups of finite irreducible Coxeter
groups. We also prove a Weyl group analogue of Cartan's theorem that all maximal tori in a connected compact Lie
group are conjugate.

\section{ Introduction}

Some years ago colleagues working in the area of statistical
mechanics asked what the maximal order Abelian subgroups of the symmetric group $S_n$
looked like. Their question arose from consideration of reducible
representations constructed from tensor products of unitary
representations arising in the statistical mechanics of systems of
$n$ quantum spins. In particular they wanted to understand the
situation as $n \rightarrow \infty$. A complete classification can
be derived from general results in \cite{Kovacs} and a classification was
given in a more general setting in \cite{Dixon}. An elementary
classification was given in \cite{Burns} using Lagrange multipliers. This
method indicates that in order to maximize the product $\prod m_{i}$
of the prime powers $m_{i}$ (the Abelian invariants) subject to
the constraint $\sum m_{i} \leq n$ (because it is an Abelian subgroup of $S_n$), all or as many as
possible of the integers $m_{i}$ should be chosen equal (to $m$
say). The problem then amounts to maximizing $m^{\frac{n}{m}}$ and regarding
this as a function of a real variable having a maximum at $e$, we
would expect that the solution to the integer-valued problem (and
therefore the maximal order  of an Abelian subgroup of $S_{n}$) is
of the form $3^{k}$, since $2^{\frac{n}{2}} < 3^{\frac{n}{3}}$.
This is essentially the case (see Theorem~\ref{thm:1.2} below).
In this note we give a complete classification of the maximal order Abelian
subgroups $M$ for all finite irreducible Coxeter groups. We also determine the number of conjugacy classes of maximal order Abelian subgroups and viewing a distinguished class of these subgroups as discrete analogues of maximal tori in compact Lie groups, we obtain a Weyl group analogue
of Cartan's theorem that all maximal tori in a connected compact Lie group $G$
are conjugate, namely:
\begin{thm*} Let $M$ and $M'$ be discrete maximal tori of $W,$ then $M'= w^{-1} Mw$ for some $w \in W.$ \\
\end{thm*}
	The first author would like  to thank the Institut de Math\'{e}matiques de Marseille,
Aix-Marseille Universit\'{e}, where part of this work was carried out and
in particular Professors Oeljeklaus and Short
for their hospitality.\\
	We now recall the precise solution for $S_n$ (i.e. $W$ of type $A_{n-1}$)
and then state the general result.
\begin{thm}\label{thm:1.2}
 Let $M$ be an Abelian subgroup of maximal order in the
symmetric group $S_n$, $n \geq 2$. Then
\begin{itemize}
\item[(i)] $M \simeq \Z_{3}^{k} \, $ if $n=3k$,

\item[(ii)] $M \simeq \Z_{3}^{k} \times  \Z_{2}$ if
$n=3k+2$,

 \item[(iii)] either $M \simeq \Z_{3}^{k-1} \times \Z_{4}$ or $M \simeq \Z_{3}^{k-1} \times  \Z_{2}
\times  \Z_{2}$ if $n=3k+1$.
\end{itemize}
\end{thm}
A natural representative for $M$ in $S_n$ is generated by a collection
of disjoint $3$-cycles, plus a $2$-cycle or a $4$-cycle if
appropriate.

Since the case of the root system of type $A_r$ is settled by the
above theorem we will exclude this case from the statement of the
general result.
\begin{thm}
Let $M$ be an Abelian subgroup of maximal order in a
finite irreducible Coxeter group $W$ of rank $r,$ then:
\begin{itemize}
\item[(a)](W crystallographic)
\begin{itemize}
\item[(i)] For $W $ of type $ B_{r}$ or $C_{r} $ we have $M \simeq
\Z_{2}^{s} \times \Z_{4}^{t}  $ where $ 0 \leq s, \, \,
t$ with $s+2t=r$ and $|M|=2^{r}$. \item[(ii)]For $W $ of type
$D_{2k} $ ($r=2k$) we have $M \simeq \Z_{2}^{2k}$  and $|M|=2^{r}$.
\item[(iii)] For $W$ of type
$D_{2k+1}$ ($r=2k +1$) we have $M \simeq \Z_{2}^{s} \times \Z_{4}^{t} $
where $ 0 \leq s, \, \, t$ with $s+2t=r-1$ and $|M|=2^{r-1}$.
\item[(iv)] For $W $ of type $ E_6$, we have $M \simeq
\Z_{3}^{3}, $ and $|M|=3^{3}$.
\item[(v)]For $W $ of type $E_{r},  \, r=7,8 $ we have $M \simeq
\Z_{2}^{r}$  and $|M|=2^{r}$. \item[(vi)] For $W $ of type
$F_4$ we have $M \simeq \Z_{2} \times \Z_{3}^{2}$, and
$|M|=2 \cdot 3^{2}$.\item[(vii)] For $W$ of type $G_2$ we have $M \simeq
\Z_{2} \times \Z_{3}$, and $|M|=2 \cdot 3$.

 \end{itemize}

\item[(b)](W noncrystallographic) \begin{itemize}\item[(i)] For
$W$ of type $H_3$ we have $M \simeq \Z_{2} \times
\Z_{5}$, and $|M|=2 \cdot 5$. \item[(ii)] For $W$ of type $H_4$ we
have $M \simeq \Z_{2} \times \Z_{5}^{2}$, and $|M|=2 \cdot 5^{2}$.
\item[(iii)] For $W$ of type $I_{2}(m), \, \, m \geq 5,$ we have $M \simeq
\Z_{m} $ and $|M|= m$.
 \end{itemize}
\end{itemize}

\end{thm}

For $W$ of type $B_n$, a natural representative for $M$ is generated
by a collection of $m$ negative $2$-cycles (having order $4$) and a
negative $1$-cycle if $n$ is odd.

For $W$ of type $D_{2k+1}$, a natural representative for $M$ comes
from a subgroup of type $W(B_{2k})$ contained in $W(D_{2k+1})$.

For $W$ of type $D_{2k}$, a natural representative for $M$ is a direct
product of $k$ groups of type $W(A_1) \times W(A_1)$.

\section{Basic facts and definitions}

All basic facts and definitions used can be found in \cite{Baki} or \cite{Hump}.
Let $(W,S)$ be an irreducible finite Coxeter system of rank $r$
with $S =\{s_{\alpha_{1}}, \ldots ,s_{\alpha_{r}}\}$ its set of
simple reflections. When $W$ is a Weyl group ($W$ crystallographic) we have
an associated connected compact Lie group $G$ (with Lie algebra
$\frakg$), containing (a fixed) maximal
torus $T$ (with Lie algebra $\frakt$) so that the Weyl group
$W = N_G(T)/T$. If $\frakg^{\C} = \frakt^{\C} \oplus \bigoplus_{\alpha \in \Phi} \frakg^{\alpha}$ is the root space decomposition of the complexification
of $\frakg$ with respect to $\frakt^{\C}$ (the
complexification of $\frakt$), then a root $\alpha$ is an element of the dual spaces $\frakt^{*}$
(pure imaginary valued) or $\frakit^{*}$ (real valued).
Since $G$ is compact the Killing form is negative definite on
$\frakt$ and gives an ( $\Ad(G)$ invariant) real inner
product $\langle\;,\;\rangle$ on the real vector spaces
$\frakit$ and $\frakit^{*}$. For $w \in N_G(T), \,
H \in \frakt$ and $ \alpha \in \frakt^{*}$ we
define $w(H)=\Ad(w)H = wHw^{-1}$ and $w(\alpha)(H) =
\alpha(\Ad(w^{-1})H),$ and since $\Ad(T)$ acts trivially on
$\frakt$ we obtain (faithful) induced actions of $W$.
Choosing a fundamental Weyl chamber in $\frakit$ we can
define positive roots $\Phi^+$ and $\{\alpha_{1}, \ldots ,\,
\alpha_{r}\}$ a basis of positive simple roots whose simple
reflections generate $W$. The fundamental weights $\{ \omega_{1},
\ldots ,\, \omega_{r} \}$ are defined by the conditions that $
    \langle \omega_i,\, 2\alpha_j\rangle :=
    \langle \alpha_j,\, \alpha_j\rangle\, \delta_{ij}$
for all $i,\, j$.
We will normalize $\langle\; ,\, \;\rangle$
 so that the highest root $\tilde{\alpha}$ has length
squared equal to two.
For $\alpha, \beta \in \Re$ we define
(integers) $n(\alpha,\beta)=\frac{2\langle \alpha,\, \beta\rangle}
{\langle \alpha,\, \alpha\rangle}$.

The Dynkin diagram $D$ is
the (multi) graph with $r$ vertices (labelled by the positive
simple roots), and $ c_{ij}c_{ji}$ edges joining $\alpha_i$ to $
\alpha_j$ where $c_{ij} = n(\alpha_i,\alpha_j)$. The extended
Dynkin diagram $ \tilde{D} $ (always labelled as in \cite{Baki}) is the graph constructed from $ D$ by
adding a new vertex $\alpha_{0} =-\tilde{\alpha}$ (the affine vertex or node) and joining it
to any vertex $\alpha_{i}$ by (the old rule of) $
n(\alpha_i,\tilde{\alpha}) \cdot n(\tilde{\alpha},\alpha_i)$
edges. We then write the coefficient $n_{i}$ over the vertex
$\alpha_{i}$ and $n_{0}=1$ over $\alpha_{0}$, where
$\tilde{\alpha}= \sum_{i=1}^r n_{i}\alpha_{i}$. Deletion of any
vertex from $ \tilde{D} $ and the edges connected to it, produces
a new (typically non connected) Dynkin diagram $D_1$ (with the
same number of vertices as $D$) of a semi-simple Lie subalgebra
$\frakg_1$ of $\frakg$. The Lie algebra $\frakg_1$
is said to obtained from $\frakg$ by an elementary operation.
Of course we can perform a new elementary operation on any of the
connected components of $D_1$. Continuing  this process we obtain
a chain of subalgebras $\frakg \supseteq \frakg_1
\supseteq \ldots \supseteq \frakg_m$, each obtained from its
predecessor by an elementary operation and any semi-simple Lie subalgebra of maximal rank is
obtained by a finite number of elementary operations (see \cite{Dynkin}, \cite{Ruben}). Note that when a diagram
of type $A_n$ occurs, an elementary operation does not change the
algebra.  Among the maximal rank Lie subalgebras are those
corresponding to maximal subgroups of maximal rank in $G$ and we recall from \cite{Wolf} the following for later use:

The fundamental simplex
$$\frakD_{0} = \{ h \in
\frakit: \alpha_{i}(h) \geq 0 \, \forall \,
i,\tilde{\alpha}(h) \leq 1 \}$$
has vertices $\{v_0, v_1, \ldots,
v_r \}$ where $v_0 = 0, \, \alpha_{i}(v_{j})=
\frac{1}{n_i}\delta_{ij}$ and it has the property that
every element of $G$ (connected and centerless) is conjugate to an
element of $\exp(2\pi i\frakD_{0})$. The conjugacy classes of maximal connected
subgroups of maximal rank in $G$ are obtained from it by a theorem of Borel and
de Siebenthal which we now recall.

\begin{thm}(\cite{Borel-deSieb},\cite{Wolf}, p.278) Let $G$ be a compact centerless simple Lie group
with fundamental simplex $\frakD_{0} = \{v_0, v_1, \ldots,
v_r \}$ and let $1 \leq i \leq r$.
\begin{itemize}
\item[(i)] Suppose that $n_i =1$,  then the centralizer of the
circle group $\{\exp(2\pi i tv_{i}): t \in \R\}$ is a maximal
connected subgroup of maximal rank in $G$ with
$$\{\alpha_{1},
\ldots ,\, \alpha_{i-1},\alpha_{i+1},\ldots ,\, \alpha_{r}\}$$
as a
system of simple roots.
\item[(ii)]  Suppose that $n_i$ is a prime
$p > 1$, then the centralizer of the element $\exp(2\pi i v_{i})$
(of order $p$) is a maximal connected subgroup of maximal rank in
$G$ with
$$\{\alpha_{0}, \ldots ,\,
\alpha_{i-1},\alpha_{i+1},\ldots ,\, \alpha_{r}\}$$
as a system of simple roots.
\item[(iii)] Every
maximal connected subgroup of maximal rank in $G$ is conjugate to
one of the above groups.
\end{itemize}
\end{thm}

Finally, the $\emph{trace}$ of a finite
Abelian group  $A= \Z_{m_1} \times \Z_{m_2} \times
\ldots \times \Z_{m_k}$ is the integer
$\Tr(A)= \sum_{i=1}^{k} m_i$ (see \cite{Hoff}). \\

\section{ Proof of Theorem 1.2}
 For $W$ of each possible type we first
prove the existence of an Abelian subgroup of the required order and
isomorphism type. We then check that $W$ contains no Abelian subgroups of
larger order, or other isomorphism types of Abelian subgroups of maximal order. For the
exceptional crystallographic types (except $G_2$) and the non-crystallographic type $H_4,$
 the check involves computer calculations using the
computer package CHEVIE for GAP~\cite{chevie}.
In order to prove the existence we
observe that if $K$ is a connected subgroup of $G$ of maximal rank
(necessarily equal to that of $G$) then any maximal torus of $K$
is also a maximal torus of $G$. The Weyl group $W(K)$ of $K$ can
therefore be identified with a subgroup of $W$ the Weyl group of
$G$.

	 For $W $ of type $B_{r}$ or $C_{r},$ the elementary operation (in the extended Dynkin diagram
of  $C_r$) corresponding to deletion of the vertex connected to the $\alpha_{0}
=-\tilde{\alpha}$ vertex of $ \tilde{D} $ (the vertex  $\alpha_1$ )
gives the Dynkin diagram $D_1$ of a semisimple subalgebra $\frakg_1$
with corresponding maximal rank subgroup of $G$ of type type
$A_1 \times C_{r-1}$. Repeating this process in the
component of $D_1$ corresponding to $C_{r-1}$ we eventually obtain
a maximal rank subgroup of type $A_1 \times A_1 \times \dots
\times A_1 $ ($r$ copies) and hence a subgroup $M \simeq
\Z_{2}^{r}$ of $W$. This sequence of elementary
operations (i.e. successive deletion of the vertex connected to
$-\tilde{\alpha}$, in successive extended Dynkin diagrams) we will call the Wolf sequence (on account of its connection to Wolf spaces, see \cite{Fino-Salamon},\cite{Wolf2}), and it also produces
a maximal rank subgroup of type
$A_1 \times A_1 \times \ldots \times A_1 $ ($r$ copies)  for types $D_{2k}$ and
$E_{r}, \, r=7,8 $ and hence a subgroup $M \simeq \Z_{2}^{r}$ of $W$.

	The maximal order Abelian
subgroups of $W$ for type $B_r$ or $C_r$
containing direct factors isomorphic to $\Z_{4} $  are also realised in the extended Dynkin diagram of $C_r$
(recalling that the Weyl group for root systems of type $B_2$ or $C_2$
is isomorphic to the dihedral group of order eight) as follows: As our first elementary operation we delete
the vertex $\alpha_2$ from $ \tilde{D}$ to obtain the Dynkin diagram $D_1$ of a
semisimple algebra $\frakg_1$ with corresponding maximal rank
subgroup of $G$ of type $C_2 \times C_{n-2}$. Repeating either
this process or taking the Wolf sequence, in the component of $D_1$
corresponding to $C_{n-2}$ we can eventually obtain any  subgroup $M \simeq \Z_{4}^{t} \times \Z_{2}^{s} $
where $ 0 \leq s, \, \, t$ with $s+2t=r$.

	In the case of $W $ of type $D_{2k+1}$ we note that not all
subgroups listed arise from subgroups of maximal rank, e.g. the
maximal order Abelian subgroup $M \simeq \Z_{4} \times
\Z_{4} $ when $W$ is of type $D_5$. However the maximal order
Abelian subgroups $M \simeq \Z_{2}^{2} \times \Z_{4} $
and $M \simeq \Z_{2}^{4}$ both arise from a maximal rank
subgroup of type $A_1 \times A_1 \times A_3$. The result follows
however from the $B_{2k}$ case by folding the $D_{2k+1}$ diagram
which comes from a regular embedding (taking a torus to a torus) of Lie groups (see \cite{Bump}, p.265).

For $W $ of type $E_6, \,F_4$ or $G_2$ the elementary
operation of deletion of the vertex $\alpha_{i}$ such that $n_i
=3,$ where $\tilde{\alpha}= \sum_{i=1}^r n_{i}\alpha_{i},$  gives a
maximal rank subgroup of type $A_2 \times A_2 \times A_2$, $A_2
\times A_2$ or $A_2$ respectively, giving rise to an Abelian
subgroup $M$ of $W$ with $M \simeq \Z_{3}^{3}$, $
\Z_{3}^{2}$ or $ \Z_{3}$ respectively. Since $ -1 \in W$ for $F_4$ and $G_2$
and the centre $Z(W)\simeq \Z_{2} \simeq \langle -1 \rangle $ (\cite{Dwyer}) we can extend these
groups by $Z(W)$ in both these cases.

 We now consider the noncrystallographic cases. For $W$
of type $I_{2}(m), \, \, m \geq 5$ (the Dihedral groups) the result is clear. For
$W$ of type $H_3$ and $H_4$ the classification of their maximal
proper subroot systems in \cite{Chen-Moody} and \cite{Doug-Pfeiff-Rohrle} gives rise to Abelian
subgroups $M$ of $W$ with $M \simeq \Z_{5}$ (from a maximal
subroot system of type $I_{2}(5)$), and $ \Z_{5}^{2}$ (from
a maximal subroot system of type $I_{2}(5) \times I_{2}(5)$) in
$H_3$ and $H_4$ respectively. Extending these groups by their
centers $Z(W) \simeq \Z_{2} $ gives the required Abelian
subgroups $ M.$

We now show that the obtained lower bounds on $|M|$
are also upper bounds and that there are no other isomorphism
types of maximal order Abelian subgroups. Before doing so we
remark (see \cite{Borel-Serre}, p.134) that a subgroup of $W$ isomorphic to
$\Z_{p}^{s}$ ( $p$ a prime) admits a faithful real representation of
dimension $r$ (on $ \frakt$) and therefore $s \leq r$ if
$p=2$ and $s \leq \frac{r}{2}$ if $p \neq 2$. There is therefore
no larger order elementary Abelian $2$-group in $ W(B_{r}), \, \,
W(D_{2k}), \, \,
W(E_7), $ or $W(E_8)$, and no larger order elementary Abelian
$3$-group in $W(E_6)$ than obtained above. In ruling out other possibilities we begin with the infinite families (i.e. classical types). We embed $M$ in a symmetric group $S_{N}$ and use the fact that we must then have
$\Tr(M) \leq N$. Here the vertices of the extended Dynkin diagram  $ \tilde{D} $ with $n_i=1$
and the maximal subgroup of maximal rank $K$ corresponding to part (i) of Theorem 2.1. play a role.
This subgroup is the isotropy subgroup of an Hermitian symmetric space $H=G/K$. Taking a maximal torus
 $T$ of $G$ to lie in $K,$ the Weyl group  $W$ acts transitively and faithfully on the fixed point set $F(T,H)$
 of the action of $T$ on $H.$  This set has cardinality equal to the Euler number $\chi(H)$ of $H$ which is equal to $2r$
when $W$ is of type $B_r $ or $ D_r$ (and its elements are pairwise
antipodal on totally geodesic 2-dimensional spheres in $H$) see \cite{Sanchez}, so that $\Tr(M) \leq 2r$.  Alternatively, we can instead of $F(T,H),$ take the weights $ \{\pm
\lambda_{1}, \ldots ,\, \pm \lambda_{r}\}$ of the vector representation of $\frakg$ for the simple Lie algebras of type $C_r$ and
$D_r.$ Now for $B_r$ and $D_r$ ($r$ even) the
center $Z(W) \simeq \langle -1\rangle $ is contained in $M$ so that the orbits
$\Omega_k$ of $M$ have even cardinality and $M$ is
the direct product of its restrictions to the orbits $\Omega_k$.
Since a transitive Abelian permutation group has order equal to
its degree we can rule out elements of order three in $M$ since they must
contribute at least six to $\Tr(M)$ and $|M|,$ whereas
$\Z_{2}^{3} $ contributes six to $\Tr(M)$ and eight to $|M|$.
The argument for ruling out higher torsion elements of $M$ other than four is similar, so that
$2^{r}$ is the maximal order of an Abelian subgroup in these
cases. Again the case of $D_{2k+1}$ follows from folding to
$B_{2k}$. The remaining large order cases, are in the exceptional families and were
verified by computer calculations.
\begin{df}
\begin{itemize}
\item[(i)] The 2-rank of $W$ is equal to the integer $r_2$ such that the maximal order of an elementary Abelian $2$-subgroup of $W$ is $2^{r_2}.$
\item[(ii)]  A pair of roots $\alpha$ and $\beta \in \Phi$ are said to be strongly orthogonal (s.o.) if $\alpha + \beta$ is not a root and $\alpha - \beta$ is not a root. A subset consisting of pairwise s.o. roots will be called a s.o. set of roots.
\end{itemize}

 \end{df}
We now have the following corollary to Theorem 1.2.
\begin{cly} Let $W$ be the Weyl group of an irreducible root system,
then the 2-rank ($r_2$) of $W$ is equal to the maximal cardinality of a set of strongly orthogonal roots and $r_2=r$ if and only if $-1 \in W.$ \end{cly}

{\bf Proof }: As in the proof of Theorem 1.2 the Wolf sequence of elementary
operations corresponding to successive deletion of the vertex connected to
$-\tilde{\alpha}$ (in successive extended Dynkin diagrams) produces a maximal rank subgroup of $G$ of type
$A_1 \times A_1 \times \ldots \times A_1 $ (with $r$ copies), and corresponding maximal order elementary
2-subgroup $M \simeq \Z_{2}^{r}$ of $W$, unless we start with or encounter a diagram of type $A_s$ with $ 2 \leq s,$ in which case $r_2 < r.$ This will occur only in diagrams of type $A_r, \, D_r, \, r$ odd, or $E_6,$ namely in those cases where $-1 \notin W,$  and then $r_2 $  takes the values ${\lfloor{(r+1)/2}\rfloor}, r-1$ and $4$ respectively. That the corresponding elementary $2$-groups are of maximal order was checked by computer for $E_6,$ follows from Theorem 1.2 for $D_r, \, r$ odd and by induction for $A_r.$ That these procedures also produce a set of s.o. roots of maximal cardinality $r_2$ follows from orthogonality in the simply-laced cases and from the classification of maximal sets of strongly orthogonal roots in \cite{Agaoka-Kaneda}, p.121 and p.127 otherwise.

\begin{df} A maximal order Abelian subgroup $M$ (of a finite irreducible Coxeter group $W$) with minimal number of Abelian invariants is called a discrete maximal torus of $W$ .
 \end{df}

{\bf  Remarks and examples}:
The Weyl group of Type $A_3$ is the symmetric group $S_4$ and it already hints at the definition of a maximal torus (of $W$) . $S_4$ has three conjugacy classes of maximal order Abelian subgroups, those of $M_1 = \langle(1234) \rangle, \, \, M_2 = \langle(12), (34) \rangle$ and $M_3 = \langle (12)(34), (13)(24) \rangle.$ Whereas $M_2$ and $M_3$ are isomorphic as abstract groups they are not as permutation groups, i.e. they are not conjugate in $S_4.$ On the other hand the cycle $(1234)$ is a Coxeter element and it generates (for $W$ of any type) a  maximal Abelian subgroup of $W$ (in this case also of maximal order) and a distinguished conjugacy class. Similarly the Weyl group of Type $B_2$ has three conjugacy classes of maximal order Abelian subgroups, two of which are isomorphic to $ \Z_{2}^{2}$ and the other (isomorphic to $\Z_{4}$) corresponding to the Coxeter element. Whereas the Abelian subgroup generated by the Coxeter element (although maximal) is no longer of maximal order in higher rank, by Theorem 1.2. the above conjugacy class phenomenon persists for classical types (excluding $D_{2k}$).

We now prove an analogue (for $W$) of Cartan's theorem that all maximal tori of a compact connected Lie group $G$ are conjugate in $G.$

\begin{thm}\label{thm:3.1} Let $M$ and $M'$ be discrete maximal tori of $W,$ then $M'= w^{-1} Mw$ for some $w \in W.$ \\
\end{thm}
{\bf Proof }: For  $W$ of type $A_r$ we note that Theorem~\ref{thm:1.2} and the definition of a maximal torus $M$ of $W$ imply that
\begin{itemize}
 \item[(i)] $M \simeq \Z_{3}^{k} \, $ if $r+1=3k$,

\item[(ii)] $M \simeq \Z_{3}^{k} \times  \Z_{2}$ if
$r+1=3k+2$,
and
 \item[(iii)]  $M \simeq \Z_{3}^{k-1} \times \Z_{4}$ if $r+1=3k+1$.
\end{itemize}
Since all direct factors correspond to disjoint cycles of length $2, \, 3$ or $4,$ (with the sum of all lengths equal to $r+1$) and at most one transposition occurring, the result follows from the fact that permutations of the same cycle type are conjugate in $S_n$.

For  $W$ of type $B_r,$ viewed as all signed permutations of $\{ 1, 2, \dots, r\},$  i.e. injective maps from $\{ 1, 2, \dots, r\}$ to
$\{\pm 1, \pm 2, \dots, \pm r\},$ with either $i$ or $-i$ in the image, elements can again be expressed in cyclic form, and the above argument generalizes. Cycles either contain both $ i$ and $-i$ (called negative cycles) and are of the form $(i_1 i_2 \dots i_k -i_1  -i_2 \dots -i_k)$ or do not contain both $ i$ and $-i$ for any $i$ (called positive cycles), and they occur in pairs of the form  $(i_1 i_2 \dots i_k)( -i_1  -i_2 \dots -i_k).$ Since (as with ordinary permutations) conjugation by $w$ of a signed permutation in cyclic form sends $i$ to $w(i),$ two signed permutations are conjugate if and only if they have the same number of positive and negative cycles of every length.
We now recall that a maximal torus $M$ of $W$  (by Theorem 1.2.) is of the form:
\begin{itemize}
 \item[(i)] $M \simeq \Z_{4}^{k} \, $ if $r=2k$ and

\item[(ii)] $M \simeq \Z_{4}^{k} \times  \Z_{2}$ if
$r=2k+1$.
\end{itemize}
 When $r=2k,$ the $k$ commuting $ \Z_{4} $ factors are negative cycles  $(i_1 i_2 , -i_1  -i_2 )$  (Coxeter elements of a $B_2$ or $C_2$ system), and we have that all maximal tori are conjugate.  When $r=2k+1,$ the argument is the same because the additional $ \Z_{2}$ factor must be a negative $1$-cycle. The  case of $D_r$ ($r$ odd) is similar.

 We next consider those cases where a maximal torus is of the form  $M \simeq \Z_{2}^{r},$ i.e. $D_{2k}$ and $E_r,  \, \, r \in \{7,8\}.$ Using the fact that for these cases $ \langle -1 \rangle = Z(W)$ must be contained in $M$ it is not hard to prove that $M$ has a set of $r$ generators, none of which is a nontrivial product of commuting reflections.
These generating reflections therefore yield a maximal set of orthogonal
roots that are in fact strongly orthogonal as the root systems are
simply-laced in the cases at hand (\cite{Agaoka-Kaneda}, p.117). However by \cite{Agaoka-Kaneda}, p.119
all maximal subsets of strongly orthogonal roots are in the same Weyl group orbit for
simply-laced root systems (the number of such $W$-orbits is the number of short simple roots) and therefore the corresponding stabilising
sets of reflection generators of the discrete maximal tori are conjugate in $W$.

Similarly in the case of $E_6,$ there is a unique $W$ orbit of sets of three orthogonal roots (\cite{Carter} p.14) and therefore (as we now show) there is a unique $W$ orbit of subroot systems of type $3A_2 =A_2 + A_2 + A_2,$ and therefore all stabilizers (the $M$'s) are conjugate. Deletion of the branch node in the extended Dynkin diagram gives a subroot system of type $3A_2$ with simple roots $\{\alpha_{1},\alpha_{3}\},$ $\{\alpha_{5},\alpha_{6}\}$ and $\{\tilde{\alpha} ,\alpha_{2}\}.$ Let $\{\alpha_{1}^\prime,\alpha_{3}^\prime\},$ $\{\alpha_{5}^\prime,\alpha_{6}^\prime\}$ and $\{\tilde{\alpha}^\prime ,\alpha_{2}^\prime\}$ be the simple roots of another $3A_2$ and by  \cite{Carter} p.14 we may assume that $ w(\alpha_{1}^\prime) = {\alpha_{1}}, \, w(\alpha_{5}^\prime) = {\alpha_{5}} $ and $w(\tilde{\alpha}^\prime) = \tilde{\alpha}$ for some $w \in W.$
 We now show that $w(\alpha_{2}^\prime)
 \in \{ \alpha_2 , \,  \tilde{\alpha}-  \alpha_2 \}.$
  Since $ \langle\tilde{\alpha}, w(\alpha_{2}^\prime)\rangle  =
 \langle w(\tilde{\alpha}^\prime), w(\alpha_{2}^\prime)\rangle  = \,    \langle\tilde{\alpha}^\prime, \alpha_{2}^\prime \rangle  \neq 0,$ we have that $b_2 \neq 0,$
 where $ w(\alpha_{2}^\prime) = \sum_{i=1}^r b_{i}\alpha_{i},$ because $\tilde{\alpha} = \omega_2.$  Similarly $ w(\alpha_{2}^\prime)$ must be orthogonal to $ \alpha_1=2\omega_1-\omega_3$ and $\alpha_3=-\omega_1+2\omega_3-\omega_4$ so that $2b_1=b_3$ and $b_1+b_4=2b_3=4b_1$ and therefore $3b_1=b_4.$ As there are only two positive roots with $\alpha_4$-coefficient equal to 3, namely $\tilde{\alpha}=\omega_2$ and $s_{\alpha_{2}}(\omega_2)=-\omega_2+\omega_4 = \tilde{\alpha}-\alpha_2$ (because the next highest root $s_{\alpha_{4}}(s_{\alpha_{2}}(\omega_2)) =-\omega_4+\omega_3+\omega_5$ has $\alpha_4$-coefficient equal to 2), we have that either $b_1=1$ and $ w(\alpha_{2}^\prime)= \tilde{\alpha}-\alpha_2 $ or $b_1=0=b_4=b_3$ and therefore  $w(\alpha_{2}^\prime) = \alpha_2 $ as the only positive root of the system of type $A_2+A_1+A_2$ (as $b_4=0$) with nonzero $\alpha_2$-coefficient. An identical argument gives $ w(\alpha_{6}^\prime)
 \in \{ \alpha_6 , \,  \alpha_5 +\alpha_6 \}$ and $  w(\alpha_{3}^\prime)
 \in \{ \alpha_3 , \,  \alpha_1 +\alpha_3 \}.$
 That all maximal order Abelian subgroups are conjugate for the cases $G_2, \, \, F_4, \, \, H_3$ and $H_4,$  follows from Sylow's Second Theorem, as the groups $ M/\langle-1\rangle$  are Sylow $p$-subgroups of $W/\langle-1\rangle$  with $p=3$ or $5.$

{\bf Remark}: The definition of a discrete maximal torus as a
maximal order Abelian subgroup with minimal number of Abelian
invariants applies to any finite group.  However, in general there is more
than one conjugacy class of them, as illustrated by $Q_8$.
A computer search of groups of small order indicates that
groups with a single conjugacy class  of maximal tori
are the exception rather than the rule.

\vspace{0.3in}\begin{minipage}{3in}
\begin{tabbing}      J. M. Burns and G. Pfeiffer
\\      School of Mathematical and Statistical Sciences, \\ National University of Ireland,
Galway, \\   Ireland

\end{tabbing}\end{minipage}

\end{document}